\documentclass[11pt,reqno,a4paper]{amsart}
\usepackage{amsmath}
\usepackage{amsfonts}
\usepackage{amssymb}
\usepackage{breqn}
\usepackage[latin1]{inputenc}
\usepackage{graphicx}

 \newtheorem{thm}{Theorem}[section]
 \newtheorem{cor}[thm]{Corollary}
 \newtheorem{lem}[thm]{Lemma}
 \newtheorem{prop}[thm]{Proposition}
 \theoremstyle{definition}
 \newtheorem{defn}[thm]{Definition}
 
 \theoremstyle{remark}
 \newtheorem{rem}[thm]{Remark}
 
 \numberwithin{equation}{section}

\def\N{{\mathbb N}}
\def\R{{\mathbb R}}
\def\C{{\mathbb C}}

\def\codiff{\mbox{codiff}}

\def\codiff{\mbox{codiff}}
\title[Propagation of regularity]
 {Propagation of regularity and positive definiteness: a constructive approach}
\author[J. Buescu]{J. Buescu}
\address{%
Dep. Matem\'{a}tica\\ FCUL and CMAFCIO \\ Portugal}

\email{jsbuescu@fc.ul.pt}

\thanks{The first author acknowledges partial support by  Fundação para a Ciência e
Tecnologia, UID/MAT/04561/2013. The third author was supported by CAPES\,-Brazil, grant 10884\,-13\,-0.
Part of this research was done while in a post doctoral program at the Universidade de Lisboa.}
\author{A. C. Paixão}
\address{Área Departamental de Matemática\\ ISEL\\ Portugal }
\email{apaixao@dem.isel.ipl.pt}

\author{C. P. Oliveira}
\address{UNIFEI \\ Instituto de Matemática e Computação\\ Brazil}
\email{oliveira@unifei.edu.br}



\date{2016-10-19}

\begin{document}

\begin{abstract}
We show that, for positive definite kernels, if specific forms of regularity (continuity, $\mathcal{S}_n$-differentiability or
holomorphy) hold locally on the diagonal, then they must hold globally on the whole domain of positive-definiteness.
This local-to-global propagation of regularity is constructively shown to be 
a consequence of the algebraic structure induced by the non-negativity of the associated bilinear forms
up to order 5. Consequences of these results for topological groups and for positive definite and exponentially 
convex functions are explored.
\end{abstract}

\maketitle

%
%
%
%
%
%
\section{Positive definite kernels and functions}
\label{sec_PDKFs}

\subsection{Introduction}
The phenomenon of {\em regularity propagation\/} for positive definite functions has recently attracted renewed attention. 
If an appropriate form of regularity (continuity, $C^{2n}$ differentiability in the real case, holomorphy in the complex case) 
holds in a neighborhood of the origin or the imaginary axis, then positive definiteness effectively sweeps it along horizontal lines, 
producing a local-to-global propagation of regularity to horizontal tubes in the whole domain. Regularity propagation 
is known in different contexts: in $\R^n$, where it follows directly \cite{bp_diff_PDF} or by properties of Fourier transforms 
\cite{don}; for complex functions, for which it follows by Fourier-Laplace transforms \cite{sas00}, \cite{bp_CPDFS};  for  
complexifications of finite-dimensional real vector spaces \cite{GL94}; and for involutive topological groups
\cite{Y98}. Proofs of regularity propagation usually rely on an appropriate integral representation determined by Bochner's theorem.

In this paper we show that regularity propagation holds in the much more general setting of {\em positive definite kernels\/}, 
for which no analog of the Bochner integral representation exists. 
Regularity propagation is thus a deeper and more fundamental property than previously considered, since it is a direct 
algebraic-analytic consequence of the definition of positive-definiteness. More 
specifically, we show constructively that, if the appropriate kind of regularity (continuity, 
$C^{2n}$-type regularity for real-variable kernels, holomorphy for complex-variable kernels) holds in a neighborhood 
of the diagonal, then the positive-definiteness condition ``sweeps" the whole domain, 
forcing regularity to propagate throughout it.

The key property for regularity propagation is the basic positive-defi\-ni\-te\-ness condition. Indeed, our 
results follow constructively from matrix positive definiteness properties of the relevant bilinear forms. 
While it is quite possible that more abstract representations such as those associated with the concept of reproducing 
kernel Hilbert space (RKHS) as introduced by Aronszajn \cite{aro} and Krein \cite{kre} may conceivably be used to provide an alternative
derivation of our results,
our methods provide purely algebraic and analytical proofs which are constructive from first principles 
and do not rely on the use of those more sophisticated functional analytic methods. 
In fact, it is possible that replicating these results relying solely on RKHS theory may involve technical subtleties 
which may render fully rigorous proofs quite involved. For instance, in \S \ref{subsec_differential_PDK} it becomes necessary to relate 
different positive definite kernels in distinct $\mathcal{S}_n$ classes to study  $C^{2n}$ propagation. In the RKHS context 
this would imply relating different RKHS consisting of functions of different finite order of differentiability 
equipped with with unknown norms, which seems a highly non-trivial task. In the same vein, the RKHS 
approach would not be able to circumvent use of (some version of) the constructive methods 
we develop in \S \ref{sec_differentiability} 
to attack differentiability, leaving the extra problem of showing, by density or more sophisticated functional analytic 
arguments, the existence of the relevant derivatives in the RKHS.

Regularity propagation for positive definite functions is not to be confused with the classical extension problem, which has 
a long and rich history \cite{krein, JorNie, rudin1, sas94, sas13} but a quite different nature. In this problem one seeks to extend 
a positive definite function beyond its domain of definition preserving positive definiteness and, if possible, regularity. 
In regularity propagation one starts with a positive definite function in a fixed domain and studies how specific types of 
regularity extend from local to global sets as a consequence of positive definiteness. There are however some  
connections between the two problems, see remark \ref{rem_final_remark}.

The paper is structured as follows. In the remainder of \S \ref{sec_PDKFs} we set up the relevant basic definitions.
In \S \ref{sec_continuity} we provide necessary and sufficient conditions for propagation of continuity 
from a neighborhood of the diagonal. As an application we characterize 
a naturally generalized class of 
shift-invariant kernels on 
topological groups. In \S \ref{sec_differentiability} we characterize regularity propagation for differentiable and 
holomorphic positive definite kernels. We conclude with an application to positive definite functions by means of 
the associated kernels.

\subsection{Basic definitions and properties}
\label{subsec_basic}

\begin{defn}
\label{def_PDK}
Let $\Omega$ be a nonempty set. A function $k : \Omega \times \Omega \to \C$ is a positive definite kernel on $\Omega$ 
if, for every $n \in \N$ and all
finite collections $\{x_j \}_{j=1}^n \subset \Omega$  and $\{ \xi_j \}_{j=1}^n \subset \C$, 
\begin{equation}
\label{eq_def_PDK}
\sum_{i, j=1}^n   k(x_i, x_j) \, \xi_i \,
\overline{\xi_j} \geq 0 
\end{equation}
that is,  if every square matrix $ [k(x_i, x_j)]_{i,j = 1}^n$ is positive semidefinite.  
\end{defn}

Consideration of order $1$ and order $2$ matrices easily leads to the following basic properties of positive definite kernels.

\begin{enumerate}
\item[($P_1$)] $k(x,x) = \overline{k(x,x)} \geq 0 \ \forall x \in \Omega$;
\item[($P_2$)] $k(x,y) = \overline{k(y,x)} \ \forall (x,y) \in \Omega^2$;
\item[($P_3$)] $|k(x,y)|^2 \leq k(x,x) \, k(y,y) \ \forall (x,y) \in \Omega^2$.
\end{enumerate}

Given a nonempty set $\Omega$, we define the {\em diagonal of  \/} $\Omega^2$ as the set 
\[ \mathcal{D}(\Omega^2) = \{ (x,x) : x \in \Omega \}. \]

Properties $P_1$--$P_3$ illustrate the central rôle played by the values assumed by $k$ along the diagonal $\mathcal{D}(\Omega^2)$ on the global 
behaviour of the kernel $k$. Consideration of higher order matrices in the appropriate topological and differentiable 
contexts has been shown to imply a number of additional properties of positive definite kernels (see \cite{bp_LAA}) which further 
stress the importance of the diagonal in the control of regularity.

In this paper we explore the consequences of the properties of $3^{\mbox \scriptsize rd}$, $4^{\mbox \scriptsize th}$ and $5^{\mbox \scriptsize th}$
order bilinear forms 
resulting from definition \ref{def_PDK} for the propagation of regularity of the kernel $k$ from the diagonal $\mathcal{D}(\Omega^2)$
to $\Omega^2$. Convenient assumptions will be made about the nature of the set $\Omega$ in order to derive conclusions for the continuous, 
$C^{2n}$-differentiable and holomorphic cases.

Natural connections exist between positive definite kernels and positive definite functions, allowing the identification of the latter 
with a subclass of the former. This fact is clear from the following definition and will be used later in this paper to derive regularity 
propagation for positive definite functions as a consequence of those obtained for positive definite kernels.

\begin{defn}
\label{def_PDK_in_groups}
Let $(G,+)$ be an Abelian group and suppose $x \rightarrow x^*$ is a map in $G$. Consider any nonempty subset $\Omega$ of $G$ and define 
$S= \Omega + \Omega^* = \{ z \in G: z = x + y^*, \, x, \, y \in \Omega \} $. 
Then a function $f: S \to \C$ is said to be {\em positive definite \/} if $k(x,y) \equiv f(x+y^*)$ is a positive definite kernel 
$k: \Omega^2 \to \C$, that is, if
\[ \sum_{i, j=1}^n   f(x_i+ x_j^*) \, \xi_i \, \overline{\xi_j} \geq 0 \]
for every $n \in \N$ and all
finite collections $\{x_j \}_{j=1}^n \subset \Omega$  and $\{ \xi_j \}_{j=1}^n \subset \C$.
\end{defn}

This 
definition of a positive definite function provides a natural generalization of the notion of shift-invariant positive definite kernel on $\Omega$.
Given a nonempty set $\Omega$, it is convenient to define the diagonal of $S = \Omega + \Omega^*$ as the set 
\[ D_\Omega (S) = \{  x + x^*: x \in \Omega \}. \]

\begin{rem}
\label{rem_diagonal} 
$D_\Omega (S)$ has zero as its unique element in examples (i) through (iii) below. However, example (iv) shows that this is not, in general, the case.
\end{rem}

We list a few concrete examples of this abstract setting where the map $x \mapsto x^*$ is an involution. Recall that, in a group $(G, \circ)$, an 
involution $*: G \to G$ is a map whose second power 
is the identity and which satisfies $(g_1 \circ g_2)^* = g_2^* \circ g_1^*$ (see e.g. \cite{berg}).

\begin{enumerate}
	\item[(i)] Taking $\Omega = G$ and $x^* = -x$, we have $S = G - G = G$. Then $f: G \to \C$ is a positive definite function if $k(x,y) = f(x-y)$ is a 
	positive definite kernel on $G$.
	\item[(ii)] In the special  case where $\Omega = \R^n$ and $x^* = -x$, we obtain the standard definition of positive definite function $f: \R^n \to \C$
	associated with the positive definite kernel $k: \R^{2n} \to \C$ defined by $k(x,y) = f(x-y)$. 
	\item[(iii)] If $\Omega \subset \R^n$ and $x^* = -x$, we have $S = \Omega - \Omega$. 
	Then $f: S \to \C$ is a positive definite function if $k(x,y) = f(x-y)$ is a positive definite kernel on $\Omega$. 
	\item[(iv)] Examples (i) to (iii) admit a natural counterpart where the $*-$map is replaced with the identity $x^* = x$, $S = \Omega + \Omega^* 
	= \Omega + \Omega$ and $k(x,y) = f(x+y)$. In this special case we will refer to $k$ as a co-positive definite function \cite{bp_CPDFS}. 
	Co-positive definite functions are also known as {\em exponentially convex functions\/} in the context of example (ii).
	In the real-variable context, these functions 
were introduced by Bernstein in 1929 \cite{ber} and have  been studied by Widder \cite{wid}, Devinatz \cite{dev2}, \cite{dev3} and other authors.
\item[(v)] Suppose that $\Omega \subset \C^n$, $x^* = -\overline{x}$ and $S = \{  z = x - \overline{y}: \ x, \, y \in \Omega \}  = \codiff(\Omega)$ 
(for the general definition of codifference sets	see \cite{bp_RC_PDF}, \cite{bp_CVPDF}). Then $f: S \to \C$ is a positive definite function iff $k(x,y) = f(x- \overline{y})$ 
is a positive definite kernel on $\Omega$.
\end{enumerate}

%
%
%
%
\section{Propagation of continuity for positive definite kernels}
\label{sec_continuity}

Suppose that $\Omega $ is a topological space and that $k$ is a positive definite kernel on $\Omega$, where $\Omega^2$ is 
endowed with the product topology.  For convenience, we define $k^y : \Omega \to \C$ by $k^y(x) = k(x,y)$ for each fixed $y \in \Omega$
and $k_x(y) = k(x,y)$ for each fixed $x \in \Omega$. 

The following lemmas will be used in the proof of this and next section's main results. 

\begin{lem}
\label{lem_continuity}
Suppose that $\mathcal{K}: \Omega^2 \to \C$ is a 
separately continuous function and that,  
for each $x_0 \in \Omega$, there exists 
a neighborhood $V(x_0)$ such that $\{ \mathcal{K}_x : x \in V(x_0)\}$ is equicontinuous in $\Omega$.   
Then $\mathcal{K}$ is continuous in $\Omega^2$.
\end{lem}

\proof Fix $(x_0,y_0)$ in $\Omega^2$. 
For any $(x,y) \in \Omega^2$ we have 
\[ | \mathcal{K}(x,y) - \mathcal{K}(x_0,y_0)| \leq |\mathcal{K}(x,y) - \mathcal{K}(x,y_0)| + |\mathcal{K}(x,y_0) - \mathcal{K}(x_0,y_0)|. \]
Observe that, for any fixed $y_0$,
$\mathcal{K}^{y_0}(x)$ is by hypothesis continuous in $\Omega$ and, in particular, at $x=x_0$. Hence, for every $\delta >0$ there exists a neighborhood $V^\delta(x_0)$
such that   
\[|\mathcal{K}^{y_0}(x) - \mathcal{K}^{y_0}(x_0)| = |\mathcal{K}(x,y_0) - \mathcal{K}(x_0,y_0)| < \delta/2\]
 whenever $x \in V^\delta(x_0)$. Now, by hypothesis, if
$x \in V(x_0)$ there exists a neighborhood $W(y_0)$, independent of $x$, such that 
\[ |\mathcal{K}_x(y) - \mathcal{K}_x(y_0)| = | \mathcal{K}(x,y) - \mathcal{K}(x,y_0)| < \delta/2 \]
for every $y \in W(y_0)$.

Then, if $(x,y) \in U(x_0,y_0) = \left( V(x_0) \cap V^\delta(x_0) \right) \times W(y_0)$, it follows that 
$|\mathcal{K}(x,y) - \mathcal{K}(x_0,y_0)| < \delta$. This finishes the proof. \qed

\begin{lem}
\label{lem_3_pt_ineq}
Let $\Omega$ be a set and suppose $k: \Omega^2 \to \C$ is a positive definite kernel. 
Then for every $x_1, x_2, x_3$ in $\Omega$, we have

\begin{equation}
\label{eq_quadradinho} 
|k(x_1, x_2) - k(x_1, x_3)|^2 \leq k(x_1,x_1) \, 2 \, \Re  \left( \frac{k(x_2, x_2) + k(x_3, x_3)}{2} - k(x_2, x_ 3) \right). 
\end{equation}
\end{lem}

\proof 
Let $\{x_1, x_2, x_3 \} \subset \Omega$. Since $k$ is a positive definite kernel, we have
\begin{equation}
\label{eq_PDK_3_by_3}
\sum_{i, j=1}^3   k(x_i, x_j) \, \xi_i \, \overline{\xi_j} \geq 0 
\end{equation}
for any $\{ \xi_j \}_{j=1}^3 \subset \C$. Choosing $\xi_1 = \eta, \, \xi_2 = 1$ 
and $\xi_3 = -1$ and using properties $P_1$ and $P_2$, condition \eqref{eq_PDK_3_by_3} may be written in the form
\begin{gather}
\label{eq_3_by_3_choice}
\begin{split}
 -k(x_1,x_1) |\eta|^2  & - 2 \Re(\eta(k(x_1,x_2) - k(x_1, x_3)) \\
 & \leq 2 \,\Re \left( \frac{k(x_2, x_2) + k(x_3, x_3)}{2} - k(x_2, x_ 3) \right)
\end{split}
\end{gather}
for every $\eta \in \C$. 

If $k(x_1, x_1) = 0$ then \eqref{eq_quadradinho} holds trivially since, in this case, it follows from property $P_3$ that 
the left hand side of \eqref{eq_3_by_3_choice} is zero.

Suppose then that $k(x_1, x_1) \neq 0$. Choosing $\eta = \dfrac{- (\overline{k(x_1, x_2) - k(x_1, x_3)})}{k(x_1, x_1)}$, 
we obtain from \eqref{eq_3_by_3_choice}

\begin{multline} 
-\frac{|k(x_1, x_2) - k(x_1, x_3)|^2}{k(x_1,x_1)} + 2 \frac{|k(x_1, x_2) - k(x_1, x_3)|^2}{k(x_1,x_1)}    \\
\leq 2 \, \Re \left( \frac{k(x_2, x_2) + k(x_3, x_3)}{2} - k(x_2, x_ 3) \right)
\end{multline}
which implies \eqref{eq_quadradinho}. \qed


\begin{thm}
\label{thm_continuity}
Let  $\Omega$ be a topological space. Suppose that $k: \Omega^2 \to \C$ is a positive definite kernel whose real part  $\Re(k)$ is continuous on the 
diagonal  $\mathcal{D}(\Omega^2)$. 
Then $k$ is continuous in $\Omega^2$.
\end{thm}

\proof Let $\{ x_1, x_2, x_3 \} \subset \Omega$. Then acording to lemma \ref{lem_3_pt_ineq}, we have


\[|k(x_1, x_2) - k(x_1, x_3)|^2 \leq k(x_1,x_1) \, 2 \, \Re  \left( \frac{k(x_2, x_2) + k(x_3, x_3)}{2} - k(x_2, x_ 3) \right). \]

%
%
%

Taking limits when $x_3 \to x_2$ and using continuity of $\Re (k)$ on the diagonal we conclude that 
$k$ is separately continuous on the second variable. From property $P_2$ we then derive that $k$ is 
separately continuous in $\Omega^2$. 
Since, according to the hypothesis, the real part of $k$ is continuous on $D(\Omega^2)$, it follows that $k(x,x)$
is bounded on some neighborhood $V(x_0)$ of any $x_0 \in \Omega$. Then there exists  $M > 0$ such that
$k(x_1, x_1) = | \Re k(x_1, x_1) | \leq M$ for every $x_1 \in V(x_0)$ and we derive from \eqref{eq_quadradinho} that

\[ |k(x_1, x_2) - k(x_1, x_3)|^2 \leq 2 \, M \, \Re  \left( \frac{k(x_2, x_2) + k(x_3, x_3)}{2} - k(x_2, x_ 3) \right). \]

Since the right hand side of the equation does not depend on $x_1$ we conclude, by taking limits when $x_3 \to x_2$ and using continuity of 
$\Re(k)$ on the diagonal, that $\{ k_x : x \in V(x_0) \}$
 is equicontinuous in $\Omega$.
Since this is true for any 
$x_0 \in \Omega$, we conclude by lemma \ref{lem_continuity} that $k$ is continuous in $\Omega^2$, as stated. \qed

\begin{rem}
Although a different proof of this result may be derived in the context of the theory of reproducing kernels, our direct approach has 
the advantage of clarifying the precise  analytical rôle of each hypothesis. In fact,  this result is frequently stated and proved, in 
the RKHS approach, by requiring cumulatively separate continuity of the kernel in both variables and continuity on the diagonal (see e.g. 
Schwartz \cite{sch}, Prop. 24, or Steinwart and Christmann \cite{s&c}, Lemma 4.29). Our proof shows that the hypothesis of separate continuity is
unnecessary, since it follows from continuity on the diagonal and joint continuity implied by continuity on the diagonal alone. The fact that 
continuity on the diagonal is, in this context, the more fundamental property may be judged from a result by Lehto~\cite{Leh} showing that 
separately continuous, bounded positive definite kernels are not necessarily continuous. It is worth noting, however, that 
in the original paper of Krein \cite{kre} on RKHS this result is stated in a rigorously equivalent form to theorem \ref{thm_continuity}.

\end{rem}

%
%
%
\subsection{Application to shift-invariant kernels on groups}
\label{subsec_groups}

A first and strai\-ght\-forward consequence of definition \ref{def_PDK_in_groups} and theorem \ref{thm_continuity} may be stated as follows.

\begin{cor}
\label{cor_cont_groups}
Suppose $(G,+)$ is a topological group and that $x \to x^*$ is a continuous map. Let $\Omega$ be an open subset of $G$, 
$S = \Omega + \Omega^*$ and suppose $f: S \to \C$ is a positive definite function. Then $f$ is continuous in $S$ if and only 
if $\Re(f)$ is continuous in 
$D_\Omega(S)$.
\end{cor}

\proof Since the ``only if" part of the assertion is trivial, we concentrate on the sufficiency assertion.
Define the positive 
definite kernel $k : \Omega^2 \to \C$ by $k(x,y) = f(x+y^*)$. Observe that $\Re(k(x,y)) = \Re(f(x+y^*))$ is continuous in $\mathcal{D}(\Omega^2)$ since
$\Re(f)$ is continuous in $D_\Omega(S)$ by hypothesis. Since $k$ is a 
positive definite kernel, we conclude by theorem \ref{thm_continuity} that $k$ is continuous in 
$\Omega^2$.

For any fixed $y \in \Omega$, define $\Omega_y = \Omega + \{y^*\}$ and $l_y : \Omega_y \to \Omega$ by $l_y(z) = z-y^*$. Notice that $\Omega_y$ is an
open set and that $l_y$ is continuous for the topologies induced by $G$ on $\Omega$ and $\Omega_y$. Now observe that 
$k(z-y^*, y) = k^y \circ l_y(z) = f_{|\Omega_y} (z)$ is continuous in $\Omega_y$ since $k$ is  continuous in $\Omega^2$.
Notice, furthermore, that since $\Omega_y$ is open, any open set of $\Omega_y$ is an open set of $S = \Omega + \Omega^* = \bigcup_{y \in \Omega} \Omega_y$.
For any open set $U \subset \C$, this imples that $f^{-1} (U) = \bigcup_{y \in \Omega} f_{|\Omega_y}^{-1}(U)$ is an open set of $S$, showing that 
$f$ is continuous in $S$ and concluding the proof. \qed

The next result follows from direct application of corollary \ref{cor_cont_groups} to example (iii) in \S \ref{sec_PDKFs} and 
is a standard result, readily found in the literature, for the case where $\Omega = \R^n$.

\begin{cor}
\label{cor_continuity_of_real_PDF}
Suppose $\Omega$ is an open subset of $\R^n$, $S = \Omega - \Omega$, and let $f: S \to \C$ be a positive definite function in the 
usual $\R^n$ variable sense (that is, taking $x^* = -x$).  Then $f$ 
is continuous in $S$ if and only if $\Re(f)$ is continuous at the origin.
\end{cor}

\proof The result 
follows immediately from corollary \ref{cor_cont_groups} by observing that, in this case, $D_\Omega(S) = \{ 0 \}$. \qed

\vspace{2mm}

In a similar way, direct application of corollary \ref{cor_cont_groups} the case of example (iv) leads to the following result.

\begin{cor}
\label{cor_continuity_of_complex_PDF}
Suppose $\Omega$  is an open subset of $\C^n$, $S = \{ z = x - \overline{y}, \, x, \, y \in \Omega \} \equiv \codiff(\Omega)$, and let $f: S \to \C$ 
be a positive definite function on $S$ in the usual $\C^n$ variable sense (i.e. with $x^* = -\overline{x}$). 
Then $f$ is continuous in $S$ if and only if $\Re(f)$ is continuous on $D_\Omega(S)$.
\end{cor}

%
%
%
%
\section{Propagation of differentiability for positive definite kernels and functions}
\label{sec_differentiability}

We now focus on differentiability properties of positive definite kernels of two real or complex variables. Results will be derived 
in a constructive way from direct consideration 
of $4^{\mbox \scriptsize th}$ and $5^{\mbox \scriptsize th}$ order square matrices 
associated with the corresponding bilinear forms on definition \ref{def_PDK} 
under suitable hypotheses on the regularity of the kernel on the diagonal. The main results will be established in \S \ref{subsec_differential_PDK}. 
Consequences for real and complex variable positive definite functions are explored in \S \ref{subsec_diff_ext_PDF}. 
Notation and terminology will be fixed in the next section along with some essential but not necessarily original results. 
Proofs are nevertheless provided when they shed light on the methods used in other sections.

\subsection{Differentiation of one or two real or complex variable functions}
\label{subsec_differentiation}

A real or complex variable function $f: D \to \C$ is said to be differentiable at a point $a \in \mbox{int} (D)$ if its derivative 

\begin{equation}
\label{eq_derivative}
\frac{df}{dx} (a) = \lim_{h \to 0} \frac{f(a+h) - f(a)}{h}
\end{equation}
is a complex number. 

In the complex variable case, differentiability of $f$ at every point of an open set $\Omega$ imples analiticity in $\Omega$ and is usually referred to as holomorphy in this set (see \cite{bp_LAA} for further details). The next definition is also very useful in this context. Given a set $U$, we denote by $U^*$ its complex conjugate, that is, $U^* = \{ \overline{z} : z \in U \}$. 

\begin{defn}
\label{def_anti-holom}
Let $D \subset \C$. A function $f: D \to \C$ is said to be anti-holomorphic on an open subset $U \subset D$ if there exists a holomorphic 
function $g: U^* \to \C$ such that $f(x) = g(\overline{x})$ for all $x \in U$. We define, for every $a \in U$, 

\begin{equation}
\label{eq_conjugate_deriv}
\frac{df}{d\overline{x}}_{|x=a} = \frac{dg}{dx}_{|x=\overline{a}}.
\end{equation}
\end{defn}

As a consequence, the corresponding higher-order differential operators are given by $\frac{d^m f}{d\overline{x^m}}_{|x=a} = \frac{d^mg}{dx^m}_{|x=\overline{a}}$.

A straightforward consequence of definition \ref{def_anti-holom} is that, for $a \in \mbox{int}(D)$,

\begin{equation}
\label{eq_higher_order}
 \frac{df}{d\overline{x}}(a) = \frac{dg}{dx}(\overline{a}) = \lim_{h \to 0} \frac{g(\overline{a}+\overline{h}) - g(\overline{a})}{\overline{h}} = 
\lim_{h \to 0} \frac{f(a+h) - f(a)}{\overline{h}}.
\end{equation}

The following characterization of the existence of the derivatives defined in \eqref{eq_derivative} and \eqref{eq_higher_order} will be essential for what follows.

\begin{lem}
\label{lem_characterization}
\begin{enumerate} 
\item[(i)]  Suppose $f: D \to \C$ is a real or complex variable function and $a \in \mbox{int} (D)$. Then $\frac{df}{dx} (a) = \lim_{h \to 0} \frac{f(a+h) - f(a)}{h}$ is
a complex number if and only if

\[ \lim_{(h,l) \to (0,0)} \frac{f(a+h)}{h} + f(a) \left( \frac{1}{l} - \frac{1}{h} \right) - \frac{f(a+l)}{l} = 0. \]

\item[(ii)] 
Suppose $f: D \to \C$ is a complex variable function and $a \in \mbox{int} (D)$. Then
$\frac{df}{d\overline{x}} (a) = \lim_{h \to 0} \frac{f(a+h) - f(a)}{\overline{h}}$ is a complex number if and only if 
\[ \lim_{(h,l) \to (0,0)} \frac{f(a+h)}{\overline{h}} + f(a) \left( \frac{1}{\overline{l}} - \frac{1}{\overline{h}} \right) - \frac{f(a+l)}{\overline{l}} = 0. \]
\end{enumerate}
\end{lem}

 \proof Since the ``only if'' part of the assertions is trivial, we focus on the sufficiency of the conditions given in (i) and (ii). 
Let $\mathcal{G}(h) = \frac{f(a+h) - f(a)}{h}$ in the first case and $\mathcal{G}(h) = \frac{f(a+h) - f(a)}{\overline{h}}$ in the second. 
In either case write $\mathcal{P}(h,l) = \mathcal{G}(h) - \mathcal{G}(l)$. Both statements will be simultaneously proved if we show that 
if $\lim_{(h,l) \to (0,0)} \mathcal{P}(h,l) = 0$, then $\lim_{h \to 0} \mathcal{G}(h)$ is necessarily a complex number.

We begin by showing that $ \mathcal{G}$ must be bounded on some (real or complex, according to the relevant case) neighborhood of zero.
In order to establish a contradiction, suppose that this is not the case. Then we can find a sequence  $l_n$ convergent to zero such that 
$ |\mathcal{G}(l_{n+1}) -  \mathcal{G} (l_n)| > 1$ for every $n \in \N$. Defining $h_n = l_{n+1}$, we have that $(h_n, l_n) \to (0,0)$ (in $\R^2$ or $\C^2$, according to the case) and that  $|\mathcal{P}(h_n,l_n)| = |\mathcal{G}(h_n) - \mathcal{G}(l_n)| = 
|\mathcal{G}(l_{n+1}) - \mathcal{G}(l_n)| > 1$ for every $n \in \N$. Hence $\lim_{(h,l) \to (0,0)} \mathcal{P}(h,l)$ cannot be zero, contradicting the hypothesis.

It therefore follows  that $\mathcal{G}$ is bounded in some neighborhood of zero. Then, for any sequence $l_n \to 0$, we conclude that 
$\mathcal{G} (l_n)$ is bounded in $\C$ and therefore admits a subsequence $\mathcal{G} (l_{b_n})$ converging to $\alpha \in \C$. Since
$(l_n, l_{b_n}) \to (0,0)$ (respectively in $\R^2$ or $\C^2$) and $\lim_{(h,l) \to (0,0)} \mathcal{P}(h,l) =0$ by hypothesis, it follows that 
$\mathcal{G} (l_n) -\mathcal{G} (l_{b_n}) \to 0$ as $n \to \infty$. Writing $\mathcal{G} (l_n) = \mathcal{G} (l_n) - \mathcal{G} (l_{b_n})
+ \mathcal{G} (l_{b_n})$, it follows that $\lim_{n \to \infty} \mathcal{G} (l_n) = \lim_{n \to \infty} \mathcal{G} (l_{b_n}) = \alpha$. 
Arbitrariness of $l_n$ implies that  $\lim_{h \to 0} \mathcal{G} (h) =\alpha \in \C$, finishing the proof. \qed

\vspace{2mm}

We now turn our attention to functions of two real or complex variables. We first concentrate on this last case, where particularly 
significant phenomena take place. It is convenient to define, given $\mathcal{U} \subset \C^2$, the sesquiconjugate set 
$\mathcal{U}^* = \{ (x,y) \in \C^2 : (x, \overline{y}) \in \mathcal{U}$ \}.

\begin{defn}
\label{def_sesquiholomorphic}
Let $\mathcal{U} \subset \C^2$ be an open set and $k: \mathcal{U} \to \C$. We say that $k(u,v)$ is sesquiholomorphic in $\mathcal{U} $
if there exists a separately holomorphic function $g: \mathcal{U}^* \to \C$ such that $k(u,v) = g(u,\overline{v})$ for all $(u,v) \in 
\mathcal{U}$, i.e., if $k(u,v)$ is separately holomorphic in $u$ and anti-holomorphic in $v$.
\end{defn}

\begin{prop}
\label{prop_sesqui_derivatives}
Let $\mathcal{U} \subset \C^2$ be an open set and suppose $k(u,v)$ is sesquiholomorphic in $\mathcal{U}$. 
Then, for all $m_1, \, m_2 \in \N$, $k$ has continuous partial derivatives 
$\frac{\partial^{m_1+m_2}}{\partial \overline{v}^{m_2} \partial u^{m_1}} k(u,v)$ 
of all orders with respect to 
the variables $u$ and $\overline{v}$ and the order of differentiation is immaterial.
\end{prop}

A proof of this fact may be found, after minor adaptations, in \cite{bp_LAA}, prop. 1.4.

The following lemma will be relevant for the use of finite increments in the context of differentiation. It is convenient to define 
\begin{equation}
\label{eq_finite_diff}
\Delta_{h,l} k(u,v) = k(u+h, v+l) - k(u, v+l) - k(u+h, v) + k(u,v)
\end{equation}
whenever the right-hand side is defined for $u, \, v, \, h, \, l \in \C$.

\begin{lem}
\label{lem_mixed_partials_finite_diff}
Let $\mathcal{U} \subset \C^2$ be an open set and suppose $k(u,v)$ is sesquiholomorphic in $\mathcal{U}$. Then, for $(u,v) \in \mathcal{U}$,
\begin{equation}
\label{eq__mixed_partials_finite_diff}
\frac{\partial^{2}}{\partial u \partial \overline{v}}  k(u,v) = \frac{\partial^{2}}{\partial \overline{v} \partial u } k(u,v) = 
\lim_{(h,l) \to (0,0)} \frac{\Delta_{h,l} k(u,v) }{h \, \overline{l}}.
\end{equation}
\end{lem}

\proof We will show that,  for separately holomorphic $g(u,v)$, the identity
\begin{equation}
\label{eq__mixed_partials_finite_hol}
\frac{\partial^{2}}{\partial u \partial v}  g(u,v) = \frac{\partial^{2}}{\partial v \partial u } g(u,v) = 
\lim_{(h,l) \to (0,0)} \frac{\Delta_{h,l} g(u,v) }{h \, l}
\end{equation}
holds. This will imply that, for $g$ such that $k(u,v) = g(u, \overline{v})$, we have, by virtue of \eqref{eq_conjugate_deriv} and 
\eqref{eq_higher_order},
\[ \frac{\partial^{2}}{\partial u \partial \overline{v}}  k(u,v) = \frac{\partial^{2}}{\partial u \partial v}  g(u,\overline{v}) \]
and
\[ \frac{\partial^{2}}{\partial \overline{v} \partial u} k(u,v) = \frac{\partial^{2}}{\partial v \partial u}  g(u,\overline{v}) =
\lim_{(h,l) \to (0,0)} \frac{\Delta_{h,\overline{l}} g(u,\overline{v}) }{h \, \overline{l}} = 
\lim_{(h,l) \to (0,0)} \frac{\Delta_{h,l} k(u,v) }{h \, \overline{l}} \]
implying the statement of the lemma.

In order to establish identity \eqref{eq__mixed_partials_finite_hol}, it is convenient to recall the following version of the finite 
increment formula for holomorphic functions (see e.g. \cite{ahl}, pg. 125). 
Suppose $f$ is analytic on an open set
$U \subset \C$, $D$ is a topological disc whose closure is contained in $U$, $C = \partial D$, $z$ and $h$ are complex numbers such that
$z$ and $z+h$ are in $D$. Then 
\begin{equation}
\label{eq_finite_incr} 
\dfrac{f(z+h) - f(z)}{h} = \dfrac{1}{2 \pi i} \int_C \dfrac{f(\zeta)}{(\zeta - z)(\zeta- z - h)} \, d\zeta.
\end{equation}  

Suppose $g: \mathcal{U}^* \to \C$ is separately holomorphic and that $(u,v) \in \mathcal{U}^*$. Let $D_1$ (resp. $D_2)$ 
be a disc centered at $U$ (resp. $v$ such that  $\overline{D_1} \times \overline{D_2} \subset \mathcal{U}^*$, $C_1 = \partial D_1$, 
$C_2 = \partial D_2$. Successive applications of formula \eqref{eq_finite_incr}  yield

\begin{gather}
  \label{eq_finite_inc_ratio}
    \begin{split}
& \dfrac{ \Delta_{h,l} g(u,v)}{h \, l} \\ 
 & \qquad = \dfrac{1}{(2 \pi i)^2}   \int_{C_2} \int_{C_1} 
\dfrac{ g(\zeta_1, \zeta_2)}{(\zeta_1 - u) (\zeta_1-u-h)(\zeta_2 - v) (\zeta_2-v-l)} d\zeta_1 \, d\zeta_2.               
	   \end{split}
    \end{gather}
Continuity of $g$ and the integral representation \eqref{eq_finite_inc_ratio} imply that

\begin{gather}
  \label{eq}
    \begin{split}
 \lim_{(h,l) \to (0,0)}  \dfrac{\Delta_{h,l} g(u,v)}{h \, l}  
 & = \dfrac{1}{(2 \pi i)^{2}}   \int_{C_{2}} \int_{C_{1}} 
\dfrac{ g(\zeta_{1}, \zeta_{2})}{(\zeta_{1} - u)^2 (\zeta_{2} - v)^2}  d\zeta_1 \, d\zeta_2   = \\
 & = \frac{\partial^2}{\partial v \partial u } g(u,v)  =  \frac{\partial^2}{\partial u \partial v } g(u,v),
 \end{split}
\end{gather}
concluding the proof. \qed

It will be essential for what follows to define the following differentiability class; for further details see e.g. \cite{bp_pdmdrki}.

\begin{defn}
\label{def_Sn(U)}
Let $\mathcal{U} \subset \R^2$ be an open set. A function  $k : \mathcal{U} \to \C$ is said to be of class $\mathcal{S}_n(\mathcal{U})$ 
if, for every $m_1 = 0, 1, \ldots, n$ and $m_2 = 0, 1, \ldots, n,$ the partial derivatives 
$\dfrac{\partial^{m_1+m_2}}{\partial v^{m_2} \, \partial u^{m_1}}  k(u,v)$ are continuous in $\mathcal{U}$.
\end{defn}

The special significance of this class of kernels in the present context is due to the following result. 

\begin{thm}
\label{thmkm_is_RK}
Let $\Omega \subseteq \R$ be an open set and $k(x,y)$ be a positive definite kernel of class 
$\mathcal{S}_n(\Omega^2)$. Then, for all $0 \leq m \leq n$, 
\[ k_m(x,y) \equiv \frac{\partial^{2m}}{\partial y^{m} \partial x^m} k(x,y) \]
is a positive definite kernel of class $\mathcal{S}_{n-m}(\Omega^2)$. 
\end{thm}

For a proof in the case where $\Omega$ is an interval, see theorem 3.6 in \cite{bp_pdmdrki}; extension to general open sets is immediate. 
Generalizations to several real and complex variables exist (see e.g. \cite{bp_LAA}) but will not be necessary in this paper.

We are now ready to construct a version of lemma \ref{lem_mixed_partials_finite_diff} for functions of two real variables
using definition \eqref{eq_finite_diff} in the appropriate way.

\begin{lem}
\label{lem_S_n(U)}
Let $\mathcal{U} \subset \R^2$ be an open set and suppose that $k : \mathcal{U} \to \C$ is of class $\mathcal{S}_1(\mathcal{U})$.
Then, for $(u,v) \in \mathcal{U} $, 
\begin{equation}
\label{eq_real_finite_diff}
\frac{\partial^{2}}{\partial u \partial v}  k(u,v) = \frac{\partial^{2}}{\partial v \partial u } k(u,v) = 
\lim_{(h,l) \to (0,0)} \frac{\Delta_{h,l} k(u,v) }{h \, l}.
\end{equation}
\end{lem}

\proof Suppose that $I_1$ and $I_2$ are open intervals centered at $u$ and $v$, respectively, such that $\overline{I_1} \times \overline{I_2} \subset \mathcal{U}  $.
For any $h, \, l$ such that $u+h \in I_1$ and $v+l \in I_1$ define
\[ \phi_1(t) = k(u+th, v+l) - k(u+th, v), \]
observe that $\phi_1$ is differentiable in an open interval containing $[0,1]$ and that 
$\Delta_{h,l} k(u,v)  = \phi_1(1) - \phi_1(0)$. By the mean value theorem for functions of one real variable 
there exists $t_1 \in \  ]0, 1[$  such that 
\[ \phi_1(1) - \phi_1(0) = \phi_1^\prime (t_1) = \left( \frac{ \partial k}{\partial u}(u+t_1 h, v+l) - \frac{ \partial k}{\partial u}(u+t_1 h, v) \right) \, h. \]
A second application of the mean value theorem to the function 
\[ \phi_2(t) = \frac{ \partial k}{\partial u} (u+t_1 h, v+tl) \]
yields 
\[ \Delta_{h,l} k(u,v)  = \frac{\partial^{2}k}{\partial u \partial v}  (u+ t_1 h,v+ t_2 l) \, h \, l\]
for some $t_2 \in ]0,1[$. Continuity of $\frac{\partial^{2}k}{\partial v \partial u}$ then implies that
\[ \lim_{(h,l) \to (0,0)} \frac{\Delta_{h,l} k(u,v)}{h \, l} =  \frac{\partial^{2}k}{\partial v \partial u} (u,v). \]
The corresponding equality for the symmetric mixed derivative $\frac{\partial^{2}k}{\partial u \partial v}$ is established 
in a completely analogous way. \qed

For the purposes of the next section, it will be convenient to establish the following consequences from lemmas \ref{lem_mixed_partials_finite_diff} and  \ref{lem_S_n(U)}.

\begin{cor}
\label{cor_psi}
In the conditions of lemma 
\ref{lem_mixed_partials_finite_diff} (resp. lemma  \ref{lem_S_n(U)}) define 
\[ \psi(h,l) = \frac{\Delta_{h,l} k(u,v)}{h \, \overline{l}}  \ \ \ \left(\mbox{resp.} \ \frac{\Delta_{h,l} k(u,v)}{h \, l}\right) \]
and 
\[ \phi(h,l) = \psi(h,h) + \psi(l,l) - \psi(h,l) - \psi(l,h). \]
Then
\[ \lim_{(h,l) \to (0,0)} \phi(h,l) = 0. \]
\end{cor}

\proof Observe that, under the conditions of lemma 
\ref{lem_mixed_partials_finite_diff} (resp. \ref{lem_S_n(U)}), we have 

\begin{gather}
  \label{eq_psi}
    \begin{split}
		\frac{\partial^2}{\partial \overline{v} \partial u } k(u,v)  & = \lim_{(h,l) \to (0,0)} \psi(h,l) = \lim_{(h,l) \to (0,0)} \psi(l,h) \\
 & =\lim_{(h,l) \to (0,0)} \psi(h,h) = \lim_{(h,l) \to (0,0)} \psi(l,l).
 \end{split}
\end{gather}
\qed

\begin{rem}
\label{rem_real_valued}
Notice that, if $u=v$, property $P_2$ yields $\psi(h,l) = \overline{\psi(l,h)}$, implying, in this case, that $\phi(h,l)$ is a real-valued function.
\end{rem}

%
%
%
%
\subsection{Propagation of differentiability for positive definite kernels}
\label{subsec_differential_PDK}

We first consider the case of complex variable positive definite kernels.

\begin{def}
\label{def_holomorphic_PDK}
Let $\Omega \subset \C$ be an open set and $k: \Omega^2 \to \C$ a positive definite kernel on $\Omega$. We say that 
$k$ is a holomorphic positive definite kernel in $ \Omega$ if $k$ is sesquiholomorphic on $\Omega^2$. 
\end{def}

\begin{rem}
Holomorphic positive definite kernels in $\Omega$, also known as holomorphic reproducing kernels, are {\em not\/} 
defined as holomorphic functions on $\Omega^2$. In fact, imposing separate holomorphy on these kernels would lead to
a trivialization of the concept; see \cite{bp_LAA} for further details.
\end{rem}

\begin{rem}
As previously observed, the theory of reproducing kernels may provide an alternative approach to the issues under discussion.
For instance, results have been stated in Krein \cite{kre} for holomorphic functions that
may be adapted to address the much more significant sesquiholomorphic case. Our proof of the following results 
does not rely on any of these functional-analytic constructions.
\end{rem}

We may now state our main theorem for the complex context.

\begin{thm}[Propagation of regularity, complex context]
\label{thm_main}
Suppose 
$\Omega$ is an open subset of $\C$ and let $k: \Omega^2 \to \C$ be a positive definite kernel. If $k$ is  
sesquiholomorphic on an open set $\mathcal{U}$ containing the diagonal $\mathcal{D}(\Omega^2)$, then $k$ is a holomorphic 
reproducing kernel in $\Omega$.
\end{thm}

\proof We will show that $k$ is separately holomorphic in the first variable and anti-holomorphic in the second 
variable for $(u,v) \subset \Omega^2$.
 For this purpose it will be sufficient to see that $\dfrac{\partial k}{\partial \overline{v}}$ exists for 
$(u_0,v_0) \in \Omega^2$ since, using property $P_2$, definition \ref{def_anti-holom} and formula \eqref{eq_higher_order},
 we have: 
 
\begin{gather*}
\begin{split}
\frac{\partial k}{\partial u}(u_0, v_0) & = \lim_{h \to 0} \frac{k(u_0+h, v_0) - k(u_0,v_0)}{h} =
\lim_{h \to 0} \overline{ \left( \frac{k(v_0 ,u_0 +h) - k(v_0 ,u_0)}{\overline{h}} \right)} \\ 
& = \overline{\frac{\partial k}{\partial \overline{v}}(v_0,u_0) }. 
\end{split}
\end{gather*}

Now consider an arbitrary
collection $\{ x_j \}_{j=1}^4 \subset \Omega$ and the corresponding bilinear form in definition 
\ref{def_PDK}. Since $k$ is a positive definite kernel, it follows that, for any collection 
$\{ \xi_j \}_{j=1}^4 \subset \C$, 
\begin{equation}
\label{eq_4-form}
\sum_{i, j=1}^4   k(x_i, x_j) \, \xi_i \,
\overline{\xi_j} \geq 0. 
\end{equation}

For $(u,v) \in \Omega^2$, fix some open disc $D \subset \C$ with center $v$ such that $\overline{D} \times  \overline{D} \subset \mathcal{U}$
and $\{ u \} \times \overline{D} \subset \Omega^2$.  For any nonzero $h,\, l\, \in \C$ such that $v + h$ and $v+l$ are in $D$, define
\[ x_1 = u, \ \  x_2 = v+h, \ \ x_3 =v, \ \ x_4 = v+l. \]
Choosing  for the $\xi_i$ the values
\[ \xi_1 = \eta, \ \ \xi_2 = \frac{1}{h}, \ \  \xi_3 = \frac{1}{l} - \frac{1}{h}, \ \ \xi_4 = -\frac{1}{l}, \ \ \ \eta \in \C, \]
and observing that  by property $P_2$ we may write 
\begin{equation}
\label{eq_hermite}
k(x_i, x_j) = \overline{k(x_j, x_i)}, \ \ \ i, \ j = 1, \ldots, 4
\end{equation}
we derive from \eqref{eq_4-form} that 
\begin{equation}
\label{eq_pos_2-form}
k(u,u) \, | \eta|^2 + 2 \Re[ \eta \, \beta_0(h,l)] + \gamma(h,l) \geq 0, 
\end{equation}
where we have defined 
\begin{equation}
\label{eq_def_beta}
\beta_0(h,l) = \frac{k(u,v+h)}{\overline{h}} + k(u,v) \left( \frac{1}{\overline{l}} - \frac{1}{\overline{h}} \right)  - \frac{k(u,v+l)}{\overline{l}}
\end{equation}
and 
\begin{gather}
\label{eq_def_gamma}
\begin{split}
\gamma(h,l)  = & \frac{k(v+h,v+h)}{|h|^2} \, + k(v,v) \left| \frac{1}{l} - \frac{1}{h} \right|^2 + \frac{k(v+l,v+l)}{|l|^2} \\
		& + 2 \, \Re \left[  k(v+h, v) \frac{1}{h} \left( \frac{1}{\overline{l}} - \frac{1}{\overline{h}} \right) 
		+ k(v+h,v+l) \, \frac{1}{h} \left( - \frac{1}{\overline{l}} \right)  \right. \\
		& \left. + k(v,v+l) 
		\left( \frac{1}{l} - \frac{1}{h} \right) \left( - \frac{1}{\overline{l}} \right)  \right].
 \end{split}
 \end{gather}

Next we show that the inequality
\begin{equation}
\label{eq_beta_gamma}
 |\beta_0(h,l) |^2 \leq k(u,u) \, \gamma(h,l).
\end{equation}
 is a  consequence of \eqref{eq_pos_2-form}. 
Indeed, if $k(u,u) = 0$, this is a trivial consequence of the fact that, by property $P_3$, $\beta_0(h,l) = 0$. 
Suppose that $k(u,u) \neq 0$. Then
\eqref{eq_beta_gamma} follows from \eqref{eq_pos_2-form} by choosing $\eta = -\frac{1}{k(u,u)} \overline{\beta_0(h,l)}$.

Finally observe that we only need to show that $\lim_{(h,l) \to (0,0)} \gamma(h,l) = 0$  in order to finish the proof.
Indeed, according to \eqref{eq_beta_gamma}, the 
aforementioned condition implies that $\lim_{(h,l) \to (0,0)} \beta_0(h,l) = 0$. 
This, in turn, will lead to the conclusion that 
$\dfrac{\partial k}{ \partial \overline{v}} (u,v) $ is a complex number, as stated in the theorem by direct application of 
lemma \ref{lem_characterization} (ii) with the appropriate identification of $k(u,v)$ with $f(v)$.

It will prove convenient 
to rewrite \eqref{eq_def_gamma} in the following form: 

\begin{gather}
\label{eq_gamma}
\begin{split}
\gamma(h,l) = & \Re \left[ \frac{k(v+h,v+h)}{|h|^2} + k(v,v) \frac{-h \overline{l} -l \overline{h} + |h|^2 + |l|^2}{ |h|^2 \,|l|^2} 
+ \frac{k(v+l,v+l)}{|l|^2}   \right. \\
	 &  \left. + 2  k(v+h, v) \left( \frac{1}{h \overline{l}} - \frac{1}{|h|^2} \right)
	-2 k(v+h,v+l) \, \frac{1}{h \overline{l}}  \right.\\
	& \left. + 2 k(v,v+l) \left( \frac{1}{h \, \overline{l}} - \frac{1}{|l|^2} \right)  \right]. 
\end{split}
\end{gather}

Now recall the contents of  corollary \ref{cor_psi} and set 
\[ \psi(h,l) = \dfrac{\Delta_{h,l} k(v,v)}{h \overline{l}}, \]
\[ \phi(h,l) = \psi(h,h) + \psi(l,l)- \psi(h,l) - \psi(l,h). \]
Then the conclusions of the corollary together with the fact that $k$ is, by hypothesis, sesquiholomorphic on an open set
$\mathcal{U}$ containing $(v,v)$, imply that
\begin{equation}
\label{eq_limit_psi}
\lim_{(h,l) \to (0,0)} \phi(h,l)  = 0.
\end{equation}

A straightforward calculation now reveals, through the use of formulas \eqref{eq_hermite} 
and direct comparison with 
\eqref{eq_gamma}, 
the following identification:

\begin{equation}
\label{eq_gamma_real}
\gamma(h,l) = \phi(h,l) 
\end{equation}
(observe that this expression  is  real by remark \ref{rem_real_valued}).  From \eqref{eq_limit_psi} and \eqref{eq_gamma_real} we finally derive 
that $\lim_{(h,l) \to (0,0)} \gamma(h,l) = 0$, which finishes the proof. \qed

\vspace{2mm}

We now focus on establishing the real variable counterpart of theorem \ref{thm_main}. The constructive 
arguments used in the proof of this theorem will be adapted in order to prove the following lemmas.

\begin{lem}
\label{thm_main_real}
Suppose $\Omega$ is an open subset of $\R$ and let $k: \Omega^2 \to \C$ be a positive definite kernel. If $k$ is of class 
$\mathcal{S}_1(\mathcal{U})$ for some open subset 
 $\mathcal{U}$ containing the diagonal $\mathcal{D}(\Omega^2)$, then $k$ is separately differentiable in $\Omega^2$.
\end{lem}

\proof The proof follows very closely that of theorem \ref{thm_main}. Since by property $P_2$ and definition \ref{def_anti-holom} 
we have that $\dfrac{\partial k}{\partial u}(u_0,v_0) = \overline{\dfrac{\partial k}{\partial v}(v_0,u_0)}$ it is sufficient to show that 
$\dfrac{\partial k}{\partial v}(u_0,v_0)$ 
exists in $\C$ for all $(u_0,v_0) \in \Omega^2$.

For $(u,v) \in \Omega^2$ we fix some open interval $I$ centered at $v$ such that 
$\overline{I} \times \overline{I} \subset \mathcal{U}$ and $ \{ u \} \times \overline{I} \subset \Omega^2$ 
and consider nonzero
$h, \, l \in \R$ such that $v+h$ and $v+l$ belong to $I$. Since $k$ is a positive definite kernel, we may write definition \ref{def_PDK} as
\eqref{eq_4-form} with the following choices:

\[ x_ 1 = u, \ \ x_2 = v + h, \ \ x_3 = v, \ \ x_4 = v+l \]
and
\[ \xi_1 = \eta, \ \ \xi_2 = \frac{1}{h}, \ \ \xi_3 = \frac{1}{l}- \frac{1}{h}, \ \ \xi_4 = -\frac{1}{l}, \ \ \ \eta \in \C. \]

Now formulas \eqref{eq_pos_2-form} through 
\eqref{eq_gamma} follow in the exact same way as in the complex variable case,
with the obvious particularity that now $\overline{h} = h$ and  $\overline{l} = l$. The arguments leading to the conclusion \eqref{eq_beta_gamma}
also carry through. Thus, if $\gamma(h,l)$ may be proved to have zero limit at the origin, direct application of lemma \ref{lem_characterization} (i) 
with the identification $k(u,v) = f(v) $ yields the conclusions of the theorem as a consequence of the fact that, in that case, 
$\lim_{(h,l) \to (0,0)} \beta_0(h,l) = 0$. 
To show that $\lim_{(h,l) \to (0,0)} \gamma(h,l) = 0$ we recall once again the definitions of $\psi$ and $\phi$ in 
corollary \ref{cor_psi} and set
\[ \psi(h,l) = \frac{\Delta_{h,l} k(v,v)}{h \, l} \]
in this case. Since $k$ is, by hypothesis, of class $\mathcal{S}_1(\mathcal{U})$ on an open set $\mathcal{U}$ containing $(v,v)$, the conclusion of the 
corollary will yield formula \eqref{eq_limit_psi} as a consequence.  The proof is finished by establishing formula \eqref{eq_gamma_real} and using it 
together with \eqref{eq_limit_psi} to conclude that $\lim_{(h,l) \to (0,0)} \gamma(h,l) = 0$. \qed 

\begin{lem}
\label{cor_real_C1}
In the conditions of lemma \ref{thm_main_real}, $\dfrac{\partial k}{\partial v}(u,v)$ 
(resp. $\dfrac{\partial k}{\partial u}(u,v)) $ 
is continuous in $v \in \Omega$ 
(resp. is continuous in $u \in \Omega)$
for every fixed value of $u \in \Omega$
(for every fixed value of $v \in \Omega$).
Furthermore, for any $(u_0,v_0) \in \Omega^2$ there exists a neighborhood $U(u_0)$ (resp. $V(v_0)$)
such that the family
$\{ \frac{\partial k}{\partial v}(u,v), \ u \in U(u_0) \}$ of fixed-$u$ functions of $v$
(resp. $\{ \frac{\partial k}{\partial u}(u,v), \ v \in V(v_0) \}$ of fixed-$v$ functions of $u$)
is equicontinuous in $\Omega$. 
\end{lem}

\proof 

As in the proof of lemma \ref{thm_main_real},  observe that by property $P_2$ we have that 
$\dfrac{\partial k}{\partial u}(u_0,v_0) = \dfrac{\overline{\partial k}}{\partial v}(v_0,u_0)$ 
for all $(u_0,v_0) \in \Omega^2$ and 
therefore it suffices to prove the assertions for $\dfrac{\partial k}{\partial v}$.

For $(u,v) \in \Omega^2$ we fix some open interval $I$ centered at $v$ such that 
$\overline{I} \times \overline{I} \subset \mathcal{U}$ and $ \{ u \} \times \overline{I} \subset \Omega^2$ 
and consider nonzero
$h, l \in \R$ such that $v+h$ and $v+h+l$ belong to $I$. Since $k$ is a positive definite kernel, we may write 
definition \ref{def_PDK} as \eqref{eq_4-form} with the following choices: 

\[ x_1 = u, \ \  x_2 = v, \ \ x_3 =v+h, \ \ x_4 = v+h+l \]
and 
\[ \xi_1 = \eta, \ \  \xi_2 = \frac{1}{h}, \ \ \xi_3 =\frac{1}{l} - \frac{1}{h}, \ \ \xi_4 = -\frac{1}{l}. \]
We formally rewrite inequality \eqref{eq_pos_2-form} as
\begin{equation}
\label{eq_pos_2-form_real}
k(u,u) \, | \eta|^2 + 2 \Re[ \eta \, \beta^0(h,l)] + \gamma^0(h,l) \geq 0, 
\end{equation}
with the new definitions

\begin{equation}
\label{eq_def_beta_zero}
\beta^0(h,l) = \frac{k(u,v+h) - k(u,v)}{h}   - \frac{k(u,v+h+l) - k(u,v+h)}{l}
\end{equation}
and 

\begin{gather}
\label{eq_def_gamma_zero}
\begin{split}
\gamma^0(h,l)  = & \frac{k(v,v)}{|h|^2} \, + k(v+h,v+h) \left| \frac{1}{l} + \frac{1}{h} \right|^2 + \frac{k(v+h+l,v+h+l)}{|l|^2} \\
		& + 2 \, \Re \left[  k(v, v+h) \left(-\frac{1}{h}\right) \left( \frac{1}{l} + \frac{1}{h} \right) 
		+ k(v,v+h+l) \, \frac{1}{hl}   \right. \\
		& \left. + k(v+h,v+h+l) 
		\left( -\frac{1}{l} \right) \left(  \frac{1}{h} + \frac{1}{l} \right)  \right].
 \end{split}
 \end{gather}

The formula corresponding to \eqref{eq_beta_gamma} is obtained in the exact same way, leading to
\begin{equation}
\label{eq_beta_gamma_real}
 |\beta^0(h,l) |^2 \leq k(u,u) \, \gamma^0(h,l).
\end{equation}
Now, defining 
\begin{equation}
\label{psi_real}
 \psi^0(h,l) = \dfrac{\Delta_{h,l} k(v,v+h)}{h \, l} 
\end{equation}
and
\begin{equation}
\label{phi_real}
 \phi^0(h,l) = \psi^0(h,h) + \psi^0(l,l)- \psi^0(h,l) - \psi^0(l,h)
\end{equation}
the identification
\begin{equation}
\label{eq_gamma_zero_real}
\gamma^0(h,l) = \phi^0(h,l) 
\end{equation}
may be established by direct comparison of \eqref{phi_real} and \eqref{eq_def_gamma_zero}. On the other hand, by adapting 
the procedures of the proof of lemma \ref{lem_S_n(U)} it is possible to derive the identity
\begin{equation}
\label{eq_mixed_partial}
\frac{\partial^2 }{\partial u \partial v} k(v,v) = \lim_{(h,l) \to (0,0)} \psi^0(h,l)
\end{equation}
from the hypothesis that $k$ is of class $\mathcal{S}_1(\mathcal{U})$ on an open set containing $(v,v)$.

From \eqref{eq_mixed_partial},      \eqref{phi_real}  and \eqref{eq_gamma_zero_real} we may now conclude that 
\[ \lim_{(h,l) \to (0,0)} \gamma^0(h,l) = \lim_{(h,l) \to (0,0)} \phi^0(h,l) = 0 \]
and also that, according to \eqref{eq_beta_gamma_real},
\begin{equation}
\label{eq_beta_zero}
\lim_{(h,l) \to (0,0)} \beta^0(h,l) = 0.
\end{equation} 
Finally, observe that condition \eqref{eq_beta_zero} implies 
\begin{equation}
\label{eq_iterated limit}
\lim_{h \to 0} \left( \lim_{l \to 0}\beta^0(h,l) \right) = 0.
\end{equation} 
Since $\lim_{l \to 0}\beta^0(h,l) = \dfrac{k(u,v+h)-k(u,v)}{h} - \dfrac{\partial k}{\partial v} (u,v+h)$, it follows from \eqref{eq_iterated limit}
that 
\[ \dfrac{\partial k}{\partial v} (u,v) = \lim_{h \to 0} \dfrac{\partial k}{\partial v} (u,v+h), \]
showing that $\dfrac{\partial k}{\partial v} (u,v) $ is continuous on the second variable $v \in \Omega$ for any fixed $u \in \Omega$, as asserted. 

To 
finish the proof we finally observe that, for any $u_0 \in \Omega$, continuity of $k(u,u)$ implies that there exists a neighborhood $U(u_0)$ such that 
$k(u,u) \leq M$ for some positive $M$ and all $u \in U(u_0)$. Hence we may rewrite formula \eqref{eq_beta_gamma_real} as 
\begin{equation}
\label{eq_beta_gamma_real_2}
 |\beta^0(h,l) |^2 \leq M \, \gamma^0(h,l).
\end{equation}
for $u \in U(u_0)$. Since the right hand side of \eqref{eq_beta_gamma_real_2} does not depend on $u$, we conclude  from 
\eqref{eq_iterated limit} that the family of fixed-$u$ functions of $v$ 
$\{ \frac{\partial k}{\partial v}(u,v), \ u \in U(u_0) \}$ is equicontinuous in $\Omega$.
\qed

The two following results will be essential in the proof of lemma \ref{lemma_2nd_deriv} below. 
The first is a corollary of lemma \ref{lem_characterization}.

\begin{cor}
\label{cor_3_14}
 Let 
$\Omega$ be an open subset of $\R$ and $k: \Omega^2 \to \C$ be a complex function. 
Suppose that $\dfrac{\partial k}{\partial u}(u,v)$ exists for every $(u,v) \in \Omega^2$
and write, whenever meaningful for $\lambda, \, h, \, l \in \R$,
\[\beta_\lambda(h,l) = \frac{k(u+\lambda, v+h)}{h} + k(u+\lambda, v) \left( \frac{1}{l} - \frac{1}{h} \right) - \frac{k(u+\lambda,v+l)}{l}. \]
Then $\dfrac{\partial^2 k}{\partial v \partial u}(u,v)$ is a complex number if and only if
\[ \lim_{(h,l) \to (0,0)} \lim_{\lambda \to 0} \frac{\beta_\lambda(h,l) - \beta_0(h,l)}{\lambda} = 0. \]
\end{cor}

\proof Applying lemma \ref{lem_characterization} to $\dfrac{\partial k}{\partial u}(u,v)$ as a function of $v$, we have that 
$\dfrac{\partial^2 k}{\partial v \partial u}(u,v)$ is a complex number if and only if
\begin{equation}
\label{eq_auxiliar}
 \lim_{(h,l) \to (0,0)} \frac{\frac{\partial k}{\partial u}(u,v+h)}{h}+\dfrac{\partial k}{\partial u}(u,v) \left( \frac{1}{l} - \frac{1}{h} \right)
- \frac{\frac{\partial k}{\partial u}(u,v+l)}{l} = 0.
\end{equation}
 Writing $\dfrac{\partial k}{\partial u}(u,v) = \lim_{\lambda \to 0}  \dfrac{k(u+\lambda,v) - k(u,v)}{\lambda}$ and using the definition of
$\beta_\lambda$, we conclude that condition \eqref{eq_auxiliar} may be written in the form
\[ \lim_{(h,l) \to (0,0)} \lim_{\lambda \to 0} \frac{\beta_\lambda(u,v) - \beta_0(u,v)}{\lambda} = 0, \]
as asserted. \qed

\begin{prop}
\label{prop_megacondensation}
Let $T$ be a square matrix of order $r_1+r_2$ partitioned in the  block form
\[ T = \begin{bmatrix} 
\frac{ 
\begin{matrix} 
\left.\begin{matrix} 
\ A \ \,
\end{matrix}\ \right| 
\left. 
\begin{matrix} 
\ B \
\end{matrix} 
\right. 
\end{matrix} 
} 
{ 
\begin{matrix} 
\left.\begin{matrix} 
\ D \ 
\end{matrix}\ \right| 
\left. 
\begin{matrix} 
\ C \ 
\end{matrix} 
\right. 
\end{matrix} 
} 
\end{bmatrix} \]
where $A = [ a_{ij} ], \, B = [ b_{iq}], \, C = [ c_{pq} ], \, D= [ d_{pj} ]$ with 
$i,j = 1\, \ldots, r_1$ and $p,q= 1, \ldots, r_2$ 
and let 
$z = (z_1, \ldots, z_{r_1}) \in \C^{r_1}, \, w = (w_1, \ldots, w_{r_2}) \in \C^{r_2}$. Then, 
if $T$ is positive semidefinite, we have 
\[ |z^T \, B \, \overline{w} |^2 \leq (z^T \,  A \, \overline{z}) \ (w^T \, C \, \overline{w}). \]
\end{prop}

References and a proof of proposition \ref{prop_megacondensation} may be found in \cite{f&h} or \cite{bp_LAA}.

\begin{lem}
\label{lemma_2nd_deriv}
Suppose $\Omega$ is an open subset of $\R$ and let $k: \Omega^2 \to \C$ be a positive definite kernel. If $k$ is of class $\mathcal{S}_1(\mathcal{U})$
for some open subset $\mathcal{U}$ containing the diagonal $\mathcal{D}(\Omega^2)$, then the second order mixed partial derivatives 
$\frac{\partial^2 k}{\partial v \partial u}(u,v)$, $\frac{\partial^2 k}{\partial u \partial v}(u,v)$ exist for all
$(u,v) \in \Omega^2$.
\end{lem}

\proof We will prove the result only for the mixed derivative $\frac{\partial^2 k}{\partial v \partial u}(u,v)$; the 
corresponding result for the mixed partial in the reverse order then follows immediately either by an analogous argument 
or simply by invoking the Hermitian property $P_2$.

For $(u,v) \in \Omega^2$ we fix some open intervals $I_u$ and $I_v$ centered at $u$ and $v$ respectively, such that 
$\overline{I_u} \times \overline{I_v} \subset \mathcal{U}$ and $\overline{I_u} \times \overline{I_v} \subset \Omega^2$. Consider
nonzero $\lambda, \, h, \, l$ such that $u+\lambda \in I_u$ and $v+l, \, v + \lambda \in I_v$ and observe that, since $k$ 
is a positive definite kernel, we may write definition \ref{def_PDK} for $n=5$ with the following choices:
\[ x_1= u, \ x_2 = u+ \lambda, \ x_3 = v + h, \  x_4 = v,  \ x_5 = v+l \]
and
\[ \xi_1= -\frac{1}{\lambda}, \ \xi_2= \frac{1}{\lambda}, \ \xi_3= \frac{1}{h}, \xi_4= \frac{1}{l} - \frac{1}{h}, \ \xi_5= -\frac{1}{l}. \]

Fix $r_1 = 2$, $r_2 = 3$ and consider the order 5 square matrix $[k(x_i, x_j)]_{i,j= 1, \ldots 5}$ partitioned in the block form given in proposition
\ref{prop_megacondensation}. This will lead to the following identifications for $i,j = 1, 2$ and $p,q = 1, 2, 3$:
\[ a_{ij} = k(x_i, x_j),  \  b_{iq} = k(x_i, x_{q+2}), \ d_{pj} = k(x_{p+2}, x_j), \ c_{pq} = k(x_{p+2}, x_{q+2}), \]
\[   z_i = \xi_i, \ w_p= \xi_{p+2}. \]
We now recall formula \eqref{eq_finite_diff} and the definitions of $\beta_\lambda$ and $\gamma$ in corollary \ref{cor_3_14} and lemma 
\ref{thm_main_real}, respectively, to obtain 
\[ z^T \, A \, \overline{z} =\frac{\Delta_{\lambda \lambda} k(u,u)}{\lambda^2}, \]
\[ z^T \, B \, \overline{w} =\frac{ \beta_\lambda(h,l) - \beta_0(h,l)}{\lambda}, \]
\[ w^T \, C \, \overline{w} = \gamma(h,l). \]
 Then, according to the conclusions of proposition \ref{prop_megacondensation}, we have:
\[ \left| \frac{\beta_\lambda (h,l) -\beta_0(h,l)}{\lambda} \right|^2 \leq \frac{\Delta_{\lambda \lambda} k(u,u)}{\lambda^2} \gamma(h,l) \]
Applying limits to both sides, we obtain:
\[ \lim_{(h,l) \to (0,0)} \lim_{\lambda \to 0} \frac{\beta_\lambda(h,l) - \beta_0(h,l)}{\lambda} \leq  
\lim_{\lambda \to 0}  \frac{\Delta_{\lambda \lambda} k(u,u)}{\lambda^2} \lim_{(h,l) \to (0,0)}  \gamma(h,l).  \]
By using lemma \ref{lem_S_n(U)}, we recognize that $\lim_{\lambda \to 0}  \frac{\Delta_{\lambda \lambda} k(u,u)}{\lambda^2}$ 
must coincide with  $\dfrac{\partial^2 k}{\partial v \partial u}(u,u)$ which, according to the hypothesis, exists for every $u \in \Omega$.
Since $ \lim_{(h,l) \to (0,0)}  \gamma(h,l) = 0$ as observed in the proof of lemma \ref{thm_main_real}, we have that 
\[ \lim_{(h,l) \to (0,0)} \lim_{\lambda \to 0} \frac{\beta_\lambda(h,l) - \beta_0(h,l)}{\lambda} = 0. \]
Hence, according to the conclusions of corollary \ref{cor_3_14}, $\dfrac{\partial^2 k}{\partial v \partial u}(u,v)$ is a 
complex number, as asserted. \qed

\begin{rem}
\label{rem_prop3.16}
It may be of interest to observe that proposition \ref{prop_megacondensation} might also have been used, as in the proof above, to 
establish the relevant inequalities in the proofs of theorem or lemma \ref{thm_main} and lemmas \ref{thm_main_real} and \ref{cor_real_C1}.
\end{rem}

We are now ready to prove our main result in the real context.

\begin{thm}[Propagation of regularity, real context]
\label{thm_main_big_real}
Suppose $\Omega$ 
is an open subset of $\R$ and let $k: \Omega^2 \to \C$ be a positive definite kernel. If $k$ is of class 
$\mathcal{S}_n(\mathcal{U})$ for some open subset $\mathcal{U}$ containing the diagonal $\mathcal{D}(\Omega^2)$, then $k$ 
is of class $\mathcal{S}_n(\Omega^2)$.
\end{thm}

\proof We first concentrate on the case $n=1$. Since $k$ is of class $\mathcal{S}_1(\mathcal{U})$, we derive from lemmas 
\ref{cor_real_C1} and \ref{lemma_2nd_deriv} that $\dfrac{\partial k}{\partial u}, \, \dfrac{\partial k}{\partial v}$,
$\dfrac{\partial^2 k}{\partial v \partial u}$
 and 
$\dfrac{\partial^2 k}{\partial u \partial v}$ exist for all $(u,v) \in \Omega^2$. Continuity of these functions will now ensure 
that $k$ is in class $\mathcal{S}_1(\Omega^2)$.

We first consider  $\dfrac{\partial k}{\partial v}$. By lemma \ref{cor_psi}, for every $u_0 \in \Omega$ there exists a neighborhood $U(u_0)$
such that the family 
$\left\{ \dfrac{\partial k}{\partial v}(u,v), \, u \in U(u_0) \right\}$ of fixed-$u$ functions of $v$ is equicontinuous in $\Omega$.
On the other hand, continuity of $\dfrac{\partial k}{\partial v} (u,v)$ in the first variable is implied by the 
existence of $\dfrac{\partial^2 k}{\partial u \partial v}$ for all $u \in \Omega$. Hence $\dfrac{\partial k}{\partial v}$ is in
the conditions of lemma \ref{lem_continuity} and we conclude that it is continuous in $\Omega^2$. A similar procedure, or
the simple observation that $k$ satisfies property $P_2$, leads to the corresponding conclusion for $\dfrac{\partial k}{\partial u}$.

In order to show that $\dfrac{\partial^2 k}{\partial u \partial v}$ is continuous in $\Omega^2$, we first recall from theorem \ref{thmkm_is_RK} 
that this function is itself 
a positive definite kernel.  Since it is, by hypothesis, continuous on $\mathcal{U}$, we conclude 
from theorem \ref{thm_continuity} that it is continuous on $\Omega^2$ and coincides (by Schwarz's theorem) with $\dfrac{\partial^2 k}{\partial v \partial u}$
and is of class $\mathcal{S}_1(\Omega^2)$. 

For the induction step, suppose the statement holds for  $n-1$. In order to prove that it also holds for $n$, suppose $k$ is in class $\mathcal{S}_n(\mathcal{U})$.
Then it is of class $\mathcal{S}_{n-1}(\mathcal{U})$ and therefore of class $\mathcal{S}_{n-1}(\Omega)$. Hence
 \[ k_{n-1} \equiv \frac{\partial^{2(n-1)} \, k}{\partial u^{n-1} \partial v^{n-1}} (u,v) \]
exists for all $(u,v) \in \Omega^2$ and, according to \cite{bp_pdmdrki}, is a positive definite kernel. Writing $k_{n-1}$ in place of $k$
and repeating the arguments of the first part of the proof, we conclude that $k_{n-1}$ has continuous first and second mixed
derivatives in $\Omega^2$. Therefore $k$ is of class $\mathcal{S}_n(\Omega^2)$, and the conclusions of the theorem now follow by induction on $n$. \qed

%
%
%
%
\subsection{Propagation of differentiability for positive definite functions}
\label{subsec_diff_ext_PDF}

We now derive consequences from theorems \ref{thm_main} and \ref{thm_main_big_real} for positive and co-positive 
definite functions. In the complex variable case, we have:

\begin{thm}
\label{thm_complex}
Suppose $\Omega \subset \C$ is an open set, $S = \codiff \ \Omega = \Omega - \Omega^*$, and  $f: S \to \C$ is a positive definite function.
If $f$ is holomorphic on an open set $U$ containing $D_\Omega (S)$, 
then $f$ is holomorphic in $S$.
\end{thm}

\proof Define $k: \Omega^2 \to \C$ by $k(u,v) = f(u-\overline{v})$. Then $k$ is a positive definite kernel. 
Continuity of the mapping $s(u,v) = u-\overline{v}$ implies that $\mathcal{U} = s^{-1}(U)$ is an open subset of $\C^2$. 
Moreover, it is clear that $\mathcal{D}(\Omega^2) \subset \mathcal{U}$. Since, for every $(x,y) \in \mathcal{U} $, we have $u-\overline{v} \in U$,
the fact that $f$ is holomorphic in $U$ implies that $k$ is sesquiholomorphic  in $\mathcal{U} $, with 

\[ \dfrac{\partial k}{\partial u} (u,v) = f^\prime(u-\overline{v}), 
\ \dfrac{\partial k}{\partial \overline{v}} (u,v) = - f^\prime(u-\overline{v}). \]
 Then, according to theorem \ref{thm_main},  we conclude that 
$k$ is a holomorphic positive definite kernel. 

Now, for any $v \in \Omega$, write $z = u-\overline{v}$ for all $u \in \Omega$ and consider the mapping 
$f: \Omega - \{\overline{v} \} \to \C$ defined by 
\[ f(z) = k(z+ \overline{v}, v). \]
This mapping is clearly holomorphic in $\Omega - \{\overline{v}\}$ since $k(u,v)$ is holomorphic in the first variable  $u \in \Omega$ for any 
$v \in \Omega$. Since $v$ is arbitrary, we conclude that $f$ is holomorphic in $S = \Omega - \Omega^*$, as asserted. \qed

\begin{thm}
\label{thm_real}
Suppose 
$\Omega \subset \R$ is an open set, $S   = \Omega - \Omega$ (resp. $S   = \Omega + \Omega)$, and  $f: S \to \C$ is a positive definite function
(resp. co-positive definite function). 
If $f$ is of class $C^{2n}$ in an open set $U$ containing 
$D_{\Omega} (S) = \{ 0 \}$ (resp. $D_{\Omega} (S)$), 
then $f$ is of class $C^{2n}$ in $S$.
\end{thm}

\proof Define $k : \Omega^2 \to \C$ by $k(u,v) = f(u+v^*)$, where $v^* = -v$ (resp.  $v^* = v$). 
 Then
$k$ is a positive definite kernel. 
Continuity of the map $s(u,v) = u + v^*$
implies that $\mathcal{U} = s^{-1} (U)$ is an open subset of $\R^2$. Moreover, it is clear that $\mathcal{D}(\Omega^2) \subset \mathcal{U}$.
Since for every $(u, v) \in \mathcal{U}$ we have $u + v^* \in U$, 
the fact that $f$ is of class $C^{2n}$ in $U$ implies that $k$ is of 
class $\mathcal{S}_n(\mathcal{U})$.  Then theorem \ref{thm_main_big_real}
implies that $k$ is of class $\mathcal{S}_n(\Omega^2)$.
Let $z_0 = u_0 + v_0^*$, where $z_0 \in S$ and $(u_0, v_0) \in \Omega^2$ and suppose $I$ is an interval containing the origin 
such that $(u_0 + I) \times (v_0 + I^*) \subset \Omega^2$, where $I^* = -I$ (resp. $I^* = I$). Then
\[ J = (u_0 + I) + (v_0 + I^*)^* = z_0 + 2I \subset S \]
is an interval containing $z_0$. For any $z \in J$, write $(u,v) = (u_0 + \frac{z-z_0}{2}, v_0 + \left( \frac{z-z_0}{2} \right )^* )$
and observe that $f(z) = k(u,v)$ to conclude that $f$ is of class $C^{2n}(J)$ with 
\[ f^{(i+j)}(z) = \left( \dfrac{1}{2} \right)^i \left( \dfrac{1^*}{2} \right)^j  \ \dfrac{\partial^{i+j} k}{\partial u^i \partial v^j} (u,v),
\ \  i, j = 1, \ldots, n. \]
Since $z_0$ is arbitrary, it follows that $f$ is of class $C^{2n}$ in $S$, as asserted. \qed

\begin{rem}
\label{rem_final_remark}
Theorem \ref{thm_complex} is particularly significant: its results may be obtained via Fourier-Laplace integral 
transforms \cite{sas00}, \cite{sas13}, \cite{bp_CPDFS} in horizontal strips of the complex plane {\em but not\/} 
in more general codifference sets. The reason our result is more general is that it does not require the existence of 
integral representations, building directly from the algebraic positive definiteness condition.

Analogously, the proof of theorem \ref{thm_real} does not require the existence of integral representations. 
For positive definite functions, in the case $\Omega = \R$ the Bochner integral representation yields equivalent results 
to our theorem \ref{thm_real} (see e.g. Donoghue \cite{don}). If $S$ does not coincide with $\R$ the situation is more delicate. 
First note that an open codifference set always contains a neighborhood $V$ of the origin. Then, given a $C^{2n} \  (n \geq 0)$  
positive definite function in $S$, 
the classical extension theorems \cite{krein}, \cite{rudin1} together with theorem \ref{thmkm_is_RK} imply that there exists a
positive definite $C^{2n}$ extension of $f_{|V}$ to $\R$. This extension is, however, in general not unique, and thus does not
allow the reconstruction of the original $f$ in $S$.

When $f$ is co-definite positive and $D_\Omega(S) = S = \Omega + \Omega$ the statements in the theorem are obviously trivial. 
However, this does not happen in general (consider, for instance, the case where $\Omega$ is the union of disjoint 
open intervals).

\end{rem}

\begin{rem}
It is worth observing that regularity propagation does {\em not\/} occur in the odd (i.e. $C^{2n+1}$) case. Wolfe \cite{wolfe}
shows constructively that, if $k$ is a positive odd integer, there exists a characteristic function $f$ such that 
$f^{(k)}(0)$ exists but $f^{(k)}(t_m)$ does not exist for a sequence of numbers $\{ t_m \}$ such that $t_m \to 0$ as 
$m \to \infty$. Characteristic functions are Fourier transforms of probability measures on $\R$, and thus, by
Bochner's theorem, coincide with the real-variable continuous positive definite functions up to a normalization factor. 
Thus this negative regularity propagation result holds for real-variable positive definite functions. 
\end{rem}

%
%
%
%


\begin{thebibliography}{99}





\bibitem{ahl}
L. Ahlfors,
\newblock {\em Complex Analysis}, 3rd ed.
\newblock McGraw-Hill, 1979.


\bibitem{aro}
N. Aronszajn, 
\newblock Theory of reproducing kernels. 
\newblock Trans. Amer. Math. Soc. 68 (1950), 337--404.

\bibitem{berg}
C. Berg, J. Christensen, P. Ressel,
\newblock {\em Harmonic analysis on semigroups.}
\newblock Springer-Verlag, GTM {\bf 100}, New York, 1984.

\bibitem{ber}
S. Bernstein, 
\newblock Sur les fonctions absolument monotones.
\newblock Acta Math. 52 (1929), 1--66.


\bibitem{sas00}
T. Bisgaard and Z. Sasvári,
\newblock{\em Characteristic functions and moment problems.}
\newblock Nova Science Publishing, New York, 2000.

\bibitem{bp_LAA}
J. Buescu, A. Paix\~{a}o,
\newblock A linear algebraic approach to holomorphic reproducing kernels in ${\mathbb C}^n$.
\newblock  Linear Algebra Appl.  412  (2006),  no. 2-3, 270--290.

\bibitem{bp_pdmdrki}
J. Buescu, A. Paix\~{a}o,
\newblock Positive definite matrices and differentiable reproducing kernel inequalities.
\newblock J. Math. Anal. Appl. {\bf 320} (2006), 279--292.

\bibitem{bp_diff_PDF}
J. Buescu, A. Paix\~{a}o,
\newblock On differentiability and analyticity of positive definite functions. 
\newblock J. Math. Anal. Appl. 375 (2011), no. 1, 336--341. 

\bibitem{bp_RC_PDF}
J. Buescu, A. Paix\~{a}o,
\newblock Real and complex variable positive definite functions. 
\newblock S\~{a}o Paulo Journal of Mathematical Sciences {\bf 6}, 2 (2012), 155--169.

\bibitem{bp_CVPDF}
J. Buescu, A. Paix\~{a}o,
\newblock Complex variable positive definite functions. 
\newblock  Complex Anal. Oper. Theory {\bf 8} (2014), no. 4, 937--954. 

\bibitem{bp_CPDFS}
J. Buescu, A. Paix\~{a}o, A. Symeonides,
\newblock Complex positive definite functions on strips.
\newblock  Complex Anal. Oper. Theory 11 (2017), no. 3, 627--649.

\bibitem{dev59}
A. Devinatz,
\newblock On the extensions of positive definite functions.
\newblock  Acta Math. 102 (1959), 109--134. 

 \bibitem{dev2}
A. Devinatz,
 \newblock Integral representations of positive definite functions.
\newblock   Trans. Amer. Math. Soc. {\bf 74} (1953), 56--77. 

\bibitem{dev3}
A. Devinatz,
\newblock Integral representations of positive definite functions II. 
\newblock  Trans. Amer. Math. Soc. {\bf 77}, (1954), 455--480.


\bibitem{don}
W. Donoghue,
\newblock {\em Distributions and Fourier transforms}.
\newblock Academic Press, New York, 1969.

\bibitem{f&h}
C. Fitzgerald, R. Horn,
\newblock On the structure of Hermitian-symmetric inequalities. 
\newblock J. Lond. Math. Soc. (2), {\bf 15} (1977), 419--430. 


\bibitem{GL94}
P. Graczyk, J. Loeb, 
\newblock Bochner and Schoenberg theorems on symmetric spaces in the complex case.
\newblock Bull. Soc. Math. France 122 (1994), no. 4, 571--590.

\bibitem{kre}
M. Krein, 
\newblock Hermitian positive kernels on homogeneous spaces I. 
\newblock Amer. Math. Soc. Transl. (2), {\bf 34} (1963), 69--108. 

\bibitem{krein}
M. Krein,
\newblock Sur le problème du prolongement des fonctions hermitiennes positives et continues. 
\newblock C. R. (Doklady) Acad. Sci. URSS (N.S.) 26 (1940). 17--22. 

\bibitem{JorNie}
P. Jorgensen, R. Niedzialomski,
\newblock Extension of positive definite functions. 
\newblock J. Math. Anal. Appl. 422 (2015), no. 1, 71--740. 

\bibitem{Leh}
O. Lehto, 
\newblock{Some remarks on the kernel functions in Hilbert spaces.}
\newblock Ann. Acad. Sci. Fenn., Ser. A I, 6 (1952), 109.


\bibitem{rudin1}
W. Rudin, 
\newblock{The extension problem for positive-definite functions.}
\newblock Illinois J. Math. 7 (1963) 532--539. 



\bibitem{sas94}
Z. Sasvári,
\newblock{\em Positive definite and definitizable functions.}
\newblock  Mathematical Topics, 2. Akademie Verlag, Berlin, 1994. 


\bibitem{sas13}
Z. Sasvári,
\newblock{\em Multivariate characteristic and correlation functions.}
\newblock De Gruyter Studies in Mathematics, 50. Walter de Gruyter \& Co., Berlin, 2013. 


\bibitem{sch}
L. Schwartz,
\newblock Sous-espaces hilbertiens d´espaces vectoriels topologiques et noyaux associés (noyaux reproduisants).
\newblock J. Analyse Math., 13 (1964), 115--256.

\bibitem{s&c}
I. Steinwart and A. Christmann, 
\newblock{ Support vector machines.}
\newblock Information Science and Statistics, Springer, New York, 2008.

\bibitem{wid} 
D. Widder,
\newblock Necessary and sufficient conditions for the representation of a function by a doubly infinite Laplace integral.
\newblock Bull. Amer. Math. Soc. {\bf 40} (1934), no. 4, 321--326.

\bibitem{wolfe}
S. Wolfe,
\newblock On derivatives of characteristic functions.
\newblock Ann. Probab. {\bf 3} (1975), 4, 737--738.


\bibitem{Y98}
E. Youssfi,
\newblock Harmonic analysis on conelike bodies and holomorphic functions on tube domains.
\newblock J. Funct. Anal. 155 (1998), no. 2, 381--435.


\end{thebibliography}
\end{document}